\documentclass[11pt]{article}
\usepackage{algorithm,algpseudocode}
\usepackage{inputenc}
\usepackage{amsmath,amssymb,epsfig,graphics,graphicx,psfrag,latexsym,amsthm}
\usepackage{subcaption}
\usepackage{fullpage,parskip,color}
\usepackage[margin=1.0in]{geometry}
\usepackage{MnSymbol,wasysym}
\usepackage[colorinlistoftodos,textsize=tiny]{todonotes}

\usepackage[normalem]{ulem}

\newtheorem{theorem}{Theorem}[section]
\newtheorem{corollary}{Corollary}[theorem]
\newtheorem{lemma}[theorem]{Lemma}

\newtheorem{definition}[theorem]{Definition}

\newcommand{\Trace}{{\hbox{Tr}}}

\newcommand{\vct}{{\hbox{vec}}}

\newcommand{\Comments}{1}
\newcommand{\mynote}[2]{\ifnum\Comments=1\textcolor{#1}{#2}\fi}
\newcommand{\mytodo}[2]{\ifnum\Comments=1%
  \todo[linecolor=#1!80!black,backgroundcolor=#1,bordercolor=#1!80!black]{#2}\fi}
\algtext*{EndWhile}
\algtext*{EndIf}
\algtext*{EndFor}

\begin{document}

\title{Metric Dimension of Hamming Graphs and Applications to Computational Biology}

\author{Lucas Laird}

\date{}

\maketitle

\begin{abstract}
Genetic sequencing has become an increasingly affordable and accessible source of genomic data in computational biology. This data is often represented as $k$-mers, i.e., strings of some fixed length $k$ with symbols chosen from a reference alphabet. In contrast, some of the most effective and well-studied machine learning algorithms require numerical representations of the data. This motivates the development of methods to embed $k$-mers into real vector spaces of low-dimension. The concept of metric dimension of the so-called Hamming graphs presents a promising way to address this issue.

A subset of vertices in a graph is said to be resolving when the distances to those vertices uniquely characterize every vertex in the graph. The metric dimension of a graph is the size of a smallest resolving subset of vertices. Unfortunately, finding the metric dimension of a general graph is a challenging problem, NP-complete in fact. Recently, however, an efficient algorithm for finding resolving sets in Hamming graphs has been proposed, which suffices to uniquely embed $k$-mers into a real vector space. Since the dimension of the embedding is the cardinality of the associated resolving set, determining whether or not a node can be removed from a resolving set while keeping it resolving is of great interest. This can be quite challenging for large graphs since only a brute-force approach is known for checking whether a set is a resolving set or not.

In this thesis, we characterize resolvability of Hamming graphs in terms of a linear system over a finite domain: a set of nodes is resolving if and only if the linear system has only a trivial solution over said domain.  Since we can represent the domain as the roots of a polynomial system, the apparatus of Gr\"obner bases comes in handy to determine, algorithmically, whether or not a set of nodes is resolving. As proof of concept, we study the resolvability of Hamming graphs associated with octapeptides i.e. proteins sequences of length eight.

This project has been directed by Prof. Manuel Lladser and co-advised by Richard C. Tillquist.
\end{abstract}

\newpage
\section{Introduction}

In recent years biological sequencing has become cheap and accessible and has brought with it an influx of genomic data. Computational biology researchers are turning to machine learning algorithms to analyze and discover patterns within these massive datasets. Analyzing these symbolic sequences is challenging however, since many machine learning algorithms require numerical representations of data. 

$k$-mers, strings of length $k$ with letters from a reference alphabet, are frequently used to represent biological sequences. Hamming graphs are an intuitive and simple method of organizing $k$-mers into a graphical format. 

We begin by formally defining the concepts of Hamming graphs and metric dimension. Hamming graphs are a specific family of connected simple graphs which are defined in terms of the Hamming distance. \\

\begin{definition}
Given two strings $s_1, s_2$ of the same length, their \textbf{Hamming distance}, denoted $d(s_1,s_2)$, is defined as the number of positions where the strings differ.\\
\end{definition}

\begin{definition}
Consider the set $S$ that is the set of all strings of length $k$ formed from a reference alphabet of size $a$; in particular,  $|S|=a^k$. The \textbf{Hamming graph} $H_{k,a}$ has vertex set $S$ such that any two vertices $x,y \in S$ are connected by an edge if and only if $d(x,y) = 1$.
\end{definition}

\begin{figure}[h]
    \centering
    \begin{subfigure}[b]{0.48\textwidth}
        \centering
        \includegraphics[width = 0.99\textwidth]{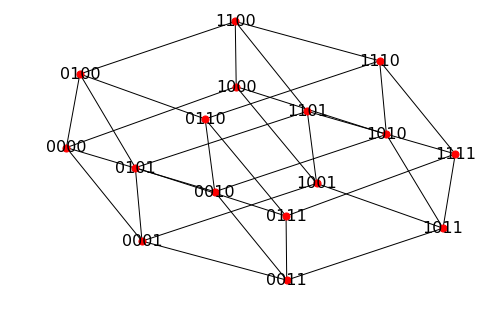}
        \caption{Hamming Graph $H_{4,2}$}
    \end{subfigure}
    \begin{subfigure}[b]{0.48\textwidth}
        \centering
        \includegraphics[width = 0.99\textwidth]{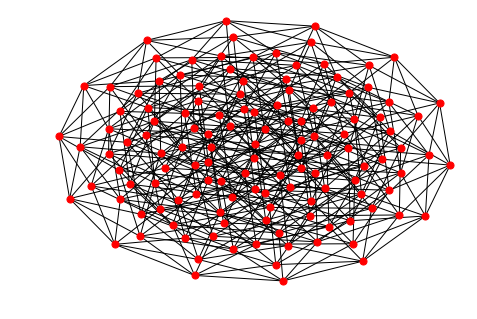}
        \caption{Hamming Graph $H_{7,2}$ (node labels omitted)}
    \end{subfigure}
    \caption{Examples of Hamming graphs}
    \label{fig:HamEx}
\end{figure}

Hamming graphs are connected and regular simple graphs where, for any two vertices, the shortest path between them is equal to their Hamming distance. Although they are simply described, they have very intricate structure and grow in size very quickly as both $a$ and $k$ increase (see Figure~\ref{fig:HamEx}).

Informally, metric dimension is the minimum number of points required to uniquely identify every point in a metric space, such as a graph, by their distances to those points. For instance, the metric dimension of the Euclidean plane is three, because the plane can be trilaterated using three non-colinear points. In the context of graphs, we have the following definition.\\

\newpage
\begin{definition}
Let $G = (V,E)$ be a connected simple graph with vertex set $V$ and edge set $E$. The distance between two vertices $d(x,y)$ is the shortest path between them. A vertex $v \in V$ \textbf{resolves} two different vertices $x,y \in V$ if $d(v,x) \neq d(v,y)$. A set $R \subseteq V$ is a \textbf{resolving set} of $G$ if every pair of vertices in $V$ is resolved by some vertex in $R$. $R$ is called a \textbf{minimal resolving set} if there is no other set smaller than $R$ that is also resolving.
\end{definition}

Note that neither resolving sets nor minimal resolving sets are unique in general.\\ 

\begin{definition}
 The \textbf{metric dimension} of a graph $G = (V,E)$, denoted $\beta(G)$, is defined as the size of a minimal resolving set.
\end{definition}

\begin{figure}[h]
    \centering
    \begin{subfigure}[b]{0.48\textwidth}
        \centering
        \includegraphics[width = 0.99\textwidth]{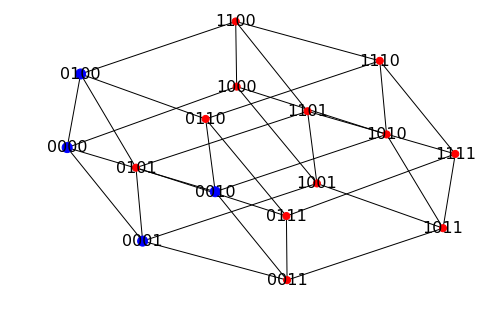}
        \caption{$H_{4,2}$ with minimal resolving set in blue}
    \end{subfigure}
    \begin{subfigure}[b]{0.48\textwidth}
        \centering
        \includegraphics[width = 0.99\textwidth]{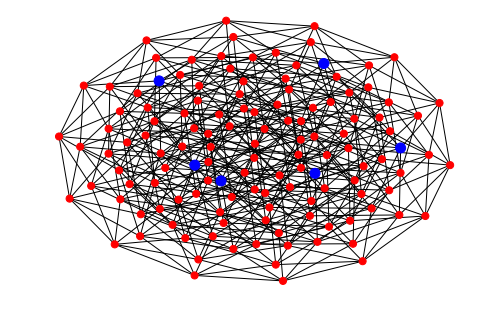}
        \caption{$H_{7,2}$ with minimal resolving set in blue}
    \end{subfigure}
    \caption{Examples of minimally resolving sets in Hamming graphs}
    \label{fig:HamResEx}
\end{figure}

For the graphs from Figure~\ref{fig:HamEx}, $\beta(H_{4,2}) = 4$ with a minimal resolving set $\{0000, 0001, 0010, 0100\}$, and $\beta(H_{7,2}) = 6$ with a minimal resolving set $\{0000000, 0000001, 0000010, 0001100, 0010100, 0100100\}$ (see Figure~\ref{fig:HamResEx}).

\section{Related Work}

Metric dimension was developed for metric spaces in 1951~\cite{Kelly_Nord:1951} and was specialized to connected simple graphs by Slater~\cite{Slater:1975} and Haray and Melter~\cite{Har_Mel:1976}. Finding the metric dimension of general connected simple graphs is a known NP-hard problem~\cite{Khu_Rag_Ros:1996, Gar_John:1990}, so most research is limited to specific families of graphs~\cite{Slater:1975, Caceres:2007} or approximation algorithms for finding nearly minimal resolving sets~\cite{Til_Lla:2018,Kratica:2007,hau_sch_vie:2012}. 

The metric dimension of Hamming graphs has been studied in many contexts, from coin-weighing problems~\cite{Guy_Now:1995} to studying optimal play strategies in the Mastermind game~\cite{Chvatal:1983}. 

Studies of metric dimension of Hamming graphs have mostly focused on hypercubes (i.e. the case with $a=2$). The exact metric dimension of these graphs are only known for strings with lengths of up to 8, and have been estimated up to length 17~ \cite{Caceres:2007,Kratica:2007}. Recently, however, Tillquist and Lladser~\cite{Til_Lla:2018} developed a computationally efficient algorithm for finding resolving sets in arbitrary Hamming graphs, which can be used to embed $k$-mers. This embedding is competitive with other state of the art graph embedding methods such as Node2Vec~\cite{grover2016node2vec} and MDS~\cite{Krz00}.
\newpage
The chief goal of this paper is to answer the question: 
\begin{center}
\textit{Given a resolving set $R$ on $H_{k,a}$, what vertices, if any, can be removed such that $R$ remains resolving?} 
\end{center}

This could be easily achieved by removing a vertex and checking if $R$ is still resolving. The only known method of checking resolvability on general Hamming graphs is a pairwise comparison of the distance vectors for every node in the graph. This brute-force approach scales as $O(n^2)$ where $n=a^k$, the number of vertices in $H_{k,a}$, making it infeasible on Hamming graphs whose vertex sets grow exponentially. Recently, it has been shown by A. F. Beardon that for hypercubes ($H_{k,2}$)~\cite{Beardon:2013}, it is possible to construct a constrained linear system for a given set $R$ such that $R$ is resolving if and only if the linear system has only a trivial solution. Unfortunately, the constraints are non-linear which makes showing that no other solutions exist quite challenging. In this paper we present a similar but distinct approach which is applicable to general Hamming graphs rather than only hypercubes. Additionally, through simplification of our system, we are able to reproduce Beardon's system on the hypercubes. Finally, we characterize the constraints of the linear system as the roots to a polynomial solution and create an efficient algorithm for checking if only the trivial solution exists.

\section{Resolvability on general Hamming graphs}

We begin by defining a useful construction.\\

\begin{definition}
The \textbf{one-hot encoding} of a vertex $v$ in $H_{k,a}$ is an $(a \times k)$ binary matrix $V$ where $V_{i,j}=1$ if $i = v[j]$, otherwise $V_{i,j}=0$.
\end{definition}

Observe that every column in a one-hot encoding has exactly one entry equal to 1 because there is only one letter represented at each position in a string.

\textbf{Example:} The one-hot encoding of 1213 in $H_{4,3}$ is the matrix $\begin{bmatrix}
1 & 0 & 1 & 0 \\
0 & 1 & 0 & 0 \\
0 & 0 & 0 & 1
\end{bmatrix}$.\\

We can use one-hot encodings to calculate the Hamming distance between two strings. In what follows, $\Trace(\cdot)$ denotes the trace of a square matrix, and $\bar A$ is the logical negation (flip) of the entries in a binary matrix $A$.\\

\begin{lemma}
For any two vertices $a$ and $b$ in $H_{k,a}$, if $A$ and $B$ denote their one-hot encodings, respectively, then the matrix $C = A^T B$ is a binary matrix such that its entries $c_{i,j} = 1$ if and only if $a[i] = b[j]$. In particular, $d(a,b) = k - \Trace(A^T B) = \Trace( A^T \bar B)$.
\label{lem:inn_prod}
\end{lemma}

Informally, this theorem states that $\Trace(A^T B)$ is counting the number of positions where $a$ and $b$ are the same. Further, logically negating $A$ results in counting the number of positions where $a$ and $b$ are not the same. A rigorous proof follows.

\noindent \emph{Proof:} Let $A_i$ and $B_i$ be the $i$-th columns of $A$ and $B$, respectively, and note that $A_i = B_j$ if and only if $a[i] = b[j]$. If $A_i = B_j$ then $A_i^T B_j = 1$ since there is exactly one $1$ in both $A_i$ and $B_j$ with the rest of the entries being $0$. If $A_i \neq B_j$ then $A_i^T B_j = 0$. Therefore, $\Trace(A^TB) = \sum_{i = 1}^k A_i^T B_i$ counts the number of positions where $a = b$, i.e. $\Trace(A^T B) = k - d(a,b)$. Similarly, if we take $\bar{B}$ to be the logical negation of $B$ then $\Trace(A^T \bar{B})$ counts the number of positions where $a \neq b$, which is by definition the Hamming distance between $a$ and $b$. \hfill$\square$

Since one-hot encodings are special types of real matrices, we note the following obvious result.\\

\begin{lemma}
Let $\vct(A)$ denote the column-wise vectorization of a matrix $A$, i.e. the vector obtained by appending the columns of $A$ to form a new column vector. For any two real valued matrices $A$ and $B$ of the same dimension, $\Trace(A^T B) = \langle \vct(A), \vct(B)\rangle$.
\label{lem:tr_vec}
\end{lemma} 

Combining Lemmas~\ref{lem:inn_prod}-\ref{lem:tr_vec}, we obtain the following general characterization of resolvability in arbitrary Hamming graphs.\\

\begin{theorem}
Let $R=\{v_1,\ldots,v_n\}$ be a subset of vertices in $H_{k,a}$ of cardinality $n$. If for each $1\le i\le n$, $V_i$ denotes the (column) vectorized one-hot encoding of $v_i$, and we define the matrix
\[A := \left(\begin{array}{c}
V_1^T \\
\vdots\\
V_n^T
\end{array}\right)\]
then $R$ is resolving if and only if there is no non-trivial solution to the linear system $Az=0$ satisfying the following constraints: if $z$ is decomposed into $k$ non-overlapping blocks of length $a$ as follows $z=$~$\big((z_1,\ldots, z_a), (z_{a+1},\ldots, z_{2a}), ... , (z_{(k-1)a+1},\ldots,z_{ka})\big)^T$ then each block is identically zero, or it has exactly one $1$ and one $(-1)$ entry and all other entries are $0$.
\label{thm:Gen_System}
\end{theorem}

\noindent \emph{Proof:} To prove the theorem, consider first two distinct vertices $x$ and $y$ in $H_{k,a}$, with vectorized one-hot encodings $X$ and $Y$, respectively. Recall that $v_i$ resolves $x$ and $y$ if and only if $d(v_i,x) \neq d(v_i,y)$, i.e. $\langle V_i,X\rangle \neq \langle V_i,Y \rangle$. Therefore, $\langle V_i, X -Y \rangle = 0$ if and only if $v_i$ does NOT resolve $x,y$. 

For $R$ to be a resolving set of $H_{k,a}$, every pair of vertices $x,y$ must be resolved by some vertex in $R$. This means R is NOT resolving if and only if there exists some distinct vertices $x,y$ such that $A(X-Y) = 0$.

Finally, since $X,Y$ are vectorized one-hot encodings, each block corresponds to a column in the one-hot encoding. There can only be one $1$ in each block for $X$ and $Y$, therefore in $z := (X-Y)$, each block can have either all $0$ entries or can have a $1$ from $X$, a $(-1)$ from $(-Y)$ and the rest of the entries are all $0$. This shows the theorem. \hfill$\square$

\textbf{Example:}
In $H_{3,2}$, consider the set of vertices $R = \{100,101,001\}$. The one-hot encodings for $R$ are: 
\[100 \leftrightarrow \begin{bmatrix} 0 & 1 & 1  \\ 1 & 0 & 0 \end{bmatrix};\quad 101 \leftrightarrow \begin{bmatrix} 0 & 1 & 0  \\ 1 & 0 & 1 \end{bmatrix};\quad 001 \leftrightarrow \begin{bmatrix} 1 & 1 & 0  \\ 0 & 0 & 1 \end{bmatrix}.\]
The matrix in Theorem~\ref{thm:Gen_System} is therefore given by 
\begin{align*}
    A := \begin{bmatrix}
    0 & 1 & 1 & 0 & 1 & 0 \\
    0 & 1 & 1 & 0 & 0 & 1 \\
    1 & 0 & 1 & 0 & 0 & 1
    \end{bmatrix}.
\end{align*}
By Theorem~\ref{thm:Gen_System}, $R$ resolves $H_{3,2}$ if and only $Az = 0$ has no non-trivial solution $z$ which satisfies the conditions from the theorem when writing $z = \big((z_1,z_2),(z_3,z_4),(z_5,z_6)\big)$. To determine whether this is the case, note the reduced row echelon form~\cite{Olver:2018} of the matrix $A$:
\begin{align*}
    \hbox{rref}(A) =  \begin{bmatrix}
    1 & 0 & 1 & 0 & 0 & 1 \\
    0 & 1 & 1 & 0 & 0 & 1 \\
    0 & 0 & 0 & 0 & 1 & -1 \\
    \end{bmatrix}.
\end{align*}

In particular, we can express the pivots $z_1$, $z_2$, and $z_5$ in terms of the free variables $z_3$, $z_4$, and $z_6$ as: $z_1 = -(z_3+z_6)$, $z_2 = -(z_3+z_6)$, and $z_5 = z_6$. Exploring now the finite number of possibilities for the free variables, we are able to determine whether there is a non-trivial solution $z$ that satisfies the constraints.

Starting with the variable $z_6$ we first recognize that $z_5,z_6$ are in the same block and $z_5 = z_6$. This means $z_5 = z_6 = 0$ since if $z_6 \neq 0$, the block constraints would require $z_5 = -z_6$. 

Then we check $z_3,z_4$ since they are in the same block and both free variables. Notice that $z_4$ is actually a useless variable since no pivots depend on it. This means we only need to consider the three possibilities for $z_3$:
\begin{itemize}
\item If $z_3 = 0$ then $z_1 = 0$, $z_2 = 0$, and $z_4 = 0$ since it is in the same block as $z_3$. This is the trivial solution $z = 0$.
\item If $z_3 = 1$ then $z_1 = -1$ and $z_2 = -1$ which is not allowed since if $z_1 = -1$, the block constraints require $z_2 = 1$.
\item If $z_3 = -1$ then $z_1 = 1$, and $z_2 = 1$ which is not allowed since if $z_1 = 1$, the block constraints require $z_2 = -1$.
\end{itemize}
As a result, there is no non-trivial solution to $Az=0$ which satisfies the constraints, so $R$ resolves $H_{3,2}$. 

This example was quite simple since no pivots depended on $z_4$, but these systems can be considerably more complex. In general, if the reduced row echelon of $A$ has $j$ free variables, then there could be up to $3^j$ possible solutions $z$ to the system $Az=0$, each of which would have to be checked for the theorem constraints. This exhaustive search could be very time consuming and difficult. Handling these constraints efficiently and algorithmically is the motivation for the next section.

\newpage
\section{Algorithmically Handling the Constraints}

Solving the linear system $Az = 0$ from Theorem~\ref{thm:Gen_System} is difficult because, while the system is linear, the constraints on the vector $z$ are not. Nevertheless, the possible $z = \big((z_1,\ldots,z_a),\ldots,(z_{(k-1)a+1},\ldots,z_{ka})\big)$ vectors can be characterized as the solutions to $p(z)=0$, with $p\in P$, where: 
\begin{equation}
P :=\left\{\begin{array}{c}
    z_1(z_1-1)(z_1+1)\\
    z_2(z_2-1)(z_2+1)\\
    \vdots\\
    z_{ka}(z_{ka}-1)(z_{ka}+1)\\
    \hline 
    \sum\limits_{i = 1}^a z_i\\
    \sum\limits_{i = a+1}^{2a} z_i\\
    \vdots \\
    \sum\limits_{i = (k-1)a+1}^{ka} z_i\\
    \hline 
    (\sum\limits_{i = 1}^a z_i^2)(2-\sum\limits_{i = 1}^a z_i^2)\\
    (\sum\limits_{i = a+1}^{2a} z_i^2)(2-\sum\limits_{i = a+1}^{2a} z_i^2)\\
    \vdots\\
    (\sum\limits_{i = (k-1)a+1}^{ka} z_i^2)(2-\sum\limits_{i = (k-1)a+1}^{ka} z_i^2)
\end{array}\right.
\label{def:P}
\end{equation}

Indeed, the set of polynomials of the form $z_i(z_i-1)(z_i+1) = 0$ enforce that $z_i \in \{-1,0,1\}$. The second set of polynomials of the form $\sum_{i = 1}^a z_i = 0$ enforce that there is a $(-1)$ for every $1$ in each block of $z$. Finally, the third set of polynomials of the form $(\sum_{i= 1}^a z_i^2)(2-\sum_{i= 1}^a z_i^2)$ ensure that there are exactly two non-zero terms or no non-zero terms in each block of $z$. These three sets of polynomials can now be combined to ensure all of the constraints on $z$ from Theorem~\ref{thm:Gen_System}.

Therefore, if we define $\{P = 0\}$ to be the set of solutions $z$ of the polynomial system, we can characterize the solution space in Theorem~\ref{thm:Gen_System} as $ker(A) \cap S$. To find $ker(A) \cap S$, we use Gaussian elimination to find $\hbox{rref}(A)$ and let $L$ be the linear equations associated with the system $rref(A)z = 0$. These linear equations are polynomials so we append them to the polynomial system $P$ since the solutions to $P \cup L$ are $ker(A) \cap \{P = 0\}$. We can show $R$ is resolving by using algebraic geometry techniques like Gr\"obner bases to check for solutions to $P \cup F$. 

\noindent\textbf{Example.} The polynomial system $P=0$ associated with the Hamming graph $H_{3,2}$ is: 
\[\left\{\begin{array}{rcl}
        z_1(z_1-1)(z_1+1) &= 0 \\
        \vdots & \\
        z_{6}(z_{6}-1)(z_{6}+1) &= 0 \\
        \hline
        z_1 + z_2 &= 0 \\
        z_3 + z_4 &= 0 \\
        z_5 + z_6 &= 0 \\
        \hline
        (z_1^2 + z_2^2 )(2-z_1^2 - z_2^2) &= 0 \\
        (z_3^2 + z_4^2)(2-z_3^2 - z_4^2) &= 0\\
        (z_5^2 + z_6^2)(2-z_5^2 - z_6^2) &= 0
    \end{array}\right.\]

\subsection{Gr{\"o}bner Bases}

Gr{\"o}bner bases are a fundamental part of algebraic geometry and are very useful to characterize solutions of polynomial systems, such as the one we encountered in the previous section. Gr{\"o}bner bases serve a purpose for polynomial ideals similar to the one orthogonal bases serve to vector spaces.

In what follows, $z=(z_1,\ldots,z_k)$ is a $k$-dimensional variable.\\

\begin{definition}
For a set of polynomials in the variable $z$, $P = \{p_1(z), \ldots , p_n(z)\}$, the associated \textit{polynomial ideal}, denoted $I(P)$, is defined as the set of all polynomials of the form $f_1 p_1 + \ldots +f_n p_n$, where $f_i(z)$ are arbitrary polynomials.
\end{definition}

Consider a homogeneous polynomial system $p_i(z) = 0$, for $i = 1,2,\ldots,n$. Define $P:=\{p_1,\ldots,p_n\}$. It follows that $z$ is a solution of the polynomial system if and only if $g(z)=0$, for all $g\in I(P)$, i.e. $z$ is a root of each polynomial in the ideal associated with the polynomials in the system. It is often more convenient to characterize the roots of the latter than those of the original system. This is done by using a Gr{\"o}bner basis of $I(P)$. Before we can define Gr{\"o}bner bases, we introduce some important concepts from~\cite{Cox_Little_OShea:2015,Cox_Little_OShea:1998}.\\

\begin{definition}[]
A \textbf{monomial} in $z$ is defined as any product of the form $z_1^{a_1}\cdots z_k^{a_k}$, with non-negative integers $a_1,\ldots,a_k$. This product is often written $z^a$ where $a=(a_1,\ldots,a_k)$.
\end{definition}

Note that all polynomials in $z$ are linear combinations of monomials.\\

\begin{definition}
A \textbf{monomial ordering} $>$ is a total ordering of the monomials such that:
\begin{enumerate}
    \item If $z^a > z^b$ then $z^a z^c > z^b z^c$ for any monomial $z^c$; and
    \item $z^a > 1$ if $z^a \neq 1$.
\end{enumerate}
\end{definition}

Some common monomial orderings are: 
\begin{itemize}
\item \textbf{Lexicographic Ordering:} $z^a > z^b$ if there is an index $i$ such that $a_j = b_j$ for all $1 \leq j < i$ but $a_i > b_i$.
\item \textbf{Graded Lexicographic Ordering: } $z^a > z^b$ if either:
\begin{enumerate}
    \item $\sum\limits_{i=1}^k a_i > \sum\limits_{i=1}^k b_i$; or
    \item $\sum\limits_{i=1}^k a_i = \sum\limits_{i=1}^k b_i$, and $z^a$ is greater than $z^b$ under the lexicographic ordering.
\end{enumerate}
\end{itemize}
The lexicographic ordering can also be reversed giving the \textbf{reversed lexicographic} and the \textbf{graded reverse lexicographic} orderings.\\

\begin{definition}
Given a polynomial $p$ and monomial ordering $>$; we can write $p = c_1 m_1 + \ldots c_n m_n$ such that $c_i$ are constants and $m_i$ are monomials in descending order, i.e. $m_1 > \ldots > m_n$. With this, we define: 
\begin{itemize}
    \item \textbf{Leading monomial}: $LM(p) := m_1$
    \item \textbf{Leading coefficient}: $LC(p) := c_1$
    \item \textbf{Leading term}: $LT(p) := c_1 m_1$
\end{itemize}
\end{definition}

\hfill\\

\begin{definition}
For a set of polynomials $P = \{p_1,\ldots,p_n\}$ and a monomial ordering $>$, every polynomial $f$ can be written through polynomial long division as: $f = q_1 p_1 + \ldots + q_n p_n + r$, where all $q_i,r$ are polynomials. Writing $f$ in this form is called \textbf{reducing} $f$ by $P$. The $q_i$ are called \textbf{quotients}, while $r$ is called the \textbf{remainder} or \textbf{reduction} of $f$ by $P$. 
\end{definition}

In what follows, we write $f\stackrel{P}{\rightarrow} r$ to mean that $r$ is the reduction of $f$ by $P$. 

The steps of multivariate polynomial long division are performed as follows. First, set $r = 0$ and $g=f$ and look for a $LM(p_i)$ that divides the $LM(g)$. If multiple $p_i$ exist we choose the one with the lowest value of $i$. With $p_i$ set $g = g - LT(g) p_i/LT(p_i)$ so that the $LT(g)$ and the $LT(p_i)$ cancel out. If no such $p_i$ exists then set $g = g-LT(g)$ and $r = r + LT(g)$. Continue this process until $g = 0$. This will produce a remainder $r$ where no monomial of $r$ is divisible by any $LM(p_i)$. 

In general, reductions are not unique because by convention we choose the $p_i$ with the lowest value of $i$. Re-ordering the polynomials in $P$ may produce a different remainder for the same polynomial $f$. This non-uniqueness is the primary motivation for Gr\"obner bases. Indeed, a Gr\"obner basis will be a set of polynomials $G$ such that for any polynomial $f$, the reduction of $f$ by $G$ is unique.\\

\begin{definition}
Given two monomials $z^a$ and $z^b$, their \textbf{least common multiple}, denoted $LCM(z^a,z^b)$, is the monomial $z^c$ with $c=(\max\{a_1,b_1\},\ldots,\max\{a_k,b_k\})$.\\
\end{definition}

\begin{definition}[]
Given two polynomials $p_1,p_2$ where $LCM(LM(p_1),LM(p_2)) = z^c$, their \textbf{S-polynomial} is defined as
\begin{align*}
    Spoly(p_1,p_2) := \frac{z^c}{LT(p_1)}p_1 - \frac{z^c}{LT(p_2)}p_2.
\end{align*}
\end{definition}

The S-polynomial is an important object for constructing Gr\"obner bases. Notice that the S-polynomial is defined so that $LT(p_1)$ and $LT(p_2)$ cancel as they do in polynomial long division. 
Also note that $Spoly(p_1,p_2) \in I(\{p_1,p_2\})$. 

With these elements in hand, Gr\"obner bases can now be defined.\\

\begin{definition} (Buchberger's Criterion.)
\label{def:buchberger_criterion}
A set of polynomials $G = \{g_1,\ldots, g_n \}$ is a \textbf{Gr{\"o}bner basis} for a polynomial ideal $I$ if and only if $I(G) = I$ and for all pairs $g_i,g_j \in G$: \[Spoly(g_1,g_2)\stackrel{G}{\rightarrow} 0.\]
A Gr{\"o}bner basis $G$ is called \textbf{reduced} if for all distinct $g_i,g_j \in G$, $LC(g_i)=1$ and no monomial of $g_j$ is divisible by $LT(g_i)$.
\end{definition}

In general, Gr\"obner bases are not unique whereas reduced Gr\"obner bases are unique under a given ordering. Gr\"obner bases also have some important properties with regards to reductions; in particular, $f\stackrel{G}{\rightarrow} 0$ if and only if $f \in I(G)$, otherwise $f\stackrel{G}{\rightarrow} r$ for some $r \notin I(G)$.

\newpage
\subsection{Computing Gr\"obner Bases}
The computation of Gr\"obner bases is an active field of research with many different classes of approaches. There are matrix reduction based algorithms, such as Faug\'ere's F4 algorithm~\cite{Faugere:1999}, as well as so-called ``signature" based algorithms, like Faug\'ere's F5 algorithm and its variants F5C and F5B~\cite{Faugere:2002,Eder:2010,Sun_Wang:2018}. The simplest and most well-studied approach is Buchberger's algorithm~\cite{Cox_Little_OShea:2015,Cox_Little_OShea:1998}. 

While significantly slower than other methods, Buchberger's algorithm (Algorithm~\ref{algo:1}) is the most easily understood and adequately demonstrates the steps involved in computing general Gr\"obner bases. 
Just like Buchberger's criterion (Definition~\ref{def:buchberger_criterion}), the algorithm is based upon S-polynomial computations. The algorithm takes in two inputs, a set of polynomials $P$ and a monomial ordering $>$. 

Algorithm~\ref{algo:1} has two bottlenecks: computing and reducing S-polynomials by the basis. Both steps rely on multivariable polynomial long division which can become very computationally intensive for high degree polynomials as well as for polynomials formed over a large variable set. Also, S-polynomials which eventually reduce to $0$ are costly to compute but contain no new information for generating the Gr\"obner basis. Avoiding the computation and reduction of S-polynomials which will reduce to $0$ is key to efficiently computing Gr\"obner bases.

\begin{algorithm}
\begin{small}
\caption{Buchberger's algorithm for computing a Gr\"obner basis}
\label{GenR_pseudo}
    \begin{algorithmic}
        \Function{buchberger}{$P$, $>$}
            \State $G$ = $P$
            \State $SP = \{Spoly(g_i,g_j)$ | $\forall i < j$,  $g_i,g_j \in G\}$
            \While{$SP$ not empty}
                \State $S$ = $SP$.pop()
                \State $r = S\stackrel{G}{\rightarrow}r$
                \If{$r \neq 0$}
                    \State Add $Spoly(r,g_i)$ to $SP$ for each $g_i \in G$ 
                    \State $G$ = $G \cup \{r\}$
                \EndIf
            \EndWhile
            \State \textbf{return} G
        \EndFunction
    \end{algorithmic}
\label{algo:1}
\end{small}
\end{algorithm}

Buchberger's algorithm may find a non-reduced Gr\"obner basis. A further basis reduction algorithm (Algorithm~\ref{algo:2}) is then applied to retrieve the unique  reduced basis for the specified ordering. Algorithm~\ref{algo:2} starts by dividing out the lead coefficient of every polynomial of $G$. Then it reduces every polynomial by the Gr\"obner basis with the polynomial removed. Since $G$ is a Gr\"obner basis, $r$ is unique and, as a result, so is the reduced Gr\"obner basis. 

\begin{algorithm}
\begin{small}
\caption{Algorithm for reducing a Gr\"obner basis}
\label{GenR_pseudo}
    \begin{algorithmic}
        \Function{reduceGr\"obner}{$G$, $>$}
            \State For each $g_i \in G$, $g_i=\frac{g_i}{LC(g_i)}$ 
            \For{$g_i \in G$}:
                \State $H=G\setminus\{g_i\}$
                \State $r=g_{i}\stackrel{H}{\rightarrow}r$
                \If{$r \neq 0$}:
                    \State $g_i = r$
                \Else
                    \State $G=G\setminus\{g_i\}$
                \EndIf
            \EndFor
            \State \textbf{return} G
        \EndFunction
    \end{algorithmic}
\label{algo:2}
\end{small}
\end{algorithm}

\subsection{Resolvability of Hamming Graphs in Terms of Gr\"obner Bases}

Gr{\"o}bner bases have many very useful properties for understanding polynomial ideals. In this work, we focus on the following important theorem.\\

\begin{theorem} (Hilbert's Weak Nullstellensatz~\cite{Cox_Little_OShea:2015,Tao:2007}.)
\label{thm:weak_null}
Let $P = \{p_1,\ldots,p_k\}$ be a set of polynomials. Then exactly one of the following statements is true:
\begin{enumerate}
    \item $P=0$ has a solution.
    \item $1\in I(P)$, i.e. there exist polynomials $f_1, \ldots, f_k$ such that $f_1 p_1+ \ldots+ f_k p_k = 1$.
\end{enumerate}
Equivalently, $\{P=0\} = \emptyset$ if and only if $\{1\}$ is the reduced Gr{\"o}bner basis of $I(P)$.
\end{theorem}

Using Hilbert's Weak Nullstellensatz, we can check for the existence of solutions to any homogeneous polynomial system; in particular, Theorem~\ref{thm:weak_null} allows us to determine whether or not the linear system $L$ together with the constraints represented as solutions to a polynomial system $P$ given by Theorem~\ref{thm:Gen_System} has solutions. This corresponds to establishing whether or not a given set of nodes in $H_{k,a}$ is resolving.
There is a slight complication to this however. $z=0$ is always a root of this linear/polynomial system. To get around this, recall that each block of $z = \big((z_1,\ldots,z_a),\ldots,(z_{(k-1)a+1},\ldots,z_{ka})\big)$ either vanishes, or contains exactly one 1 and one (-1) entry while the remaining entries vanish. As a result, if $z\ne 0$ then $\sum_{i = 1}^{ka} z_i^2=2i$ for some $1\le i\le k$. This motivates to consider the polynomials $f_i = \sum_{i = 1}^{ka} z_i^2-2i$ with $i = 1,\ldots,k$. 

If we now define $P_i = P \cup f_i$ for $i = 1,\ldots,k$, then each $P_i$ will have only non-trivial solutions---if any at all. This means that computing the reduced Gr\"obner bases $G_i$ of $Az \cup P_i$ and applying Hilbert's Weak Nullstellensatz to each $G_i$ individually can be used to check for the existence of non-trivial solutions. If any $G_i \neq \{1\}$, then there is a non-trivial solution to $Az \cup P = 0$ and $R$ is not resolving. 

\noindent\textbf{Example.} Recall from our previous example that $R = \{100,101,001\}$ resolves $H_{3,2}$. Moreover, the reduced row echelon form of the matrix $A$ given by Theorem~\ref{thm:Gen_System} is:
\begin{align*}
    \hbox{rref}(A) = \begin{bmatrix}
    1 & 0 & 1 & 0 & 0 & 1 \\
    0 & 1 & 1 & 0 & 0 & 1 \\
    0 & 0 & 0 & 0 & 1 & -1 \\
    \end{bmatrix}.
\end{align*}
and the polynomials associated with the constraints on $z$ in Theorem~\ref{thm:Gen_System} is
\[P := \left\{\begin{array}{c}
        z_1(z_1-1)(z_1+1)\\
        \vdots \\
        z_{6}(z_{6}-1)(z_{6}+1) \\
        \hline
        z_1 + z_2 \\
        z_3 + z_4 \\
        z_5 + z_6 \\
        \hline
        (z_1^2 + z_2^2 )(2-z_1^2 - z_2^2) \\
        (z_3^2 + z_4^2)(2-z_3^2 - z_4^2) \\
        (z_5^2 + z_6^2)(2-z_5^2 - z_6^2) 
\end{array}\right..\]
\newpage
Finally, we must consider the polynomials:
\[f_i = \sum\limits_{j=1}^{6} z_j^2 - 2i,\]
for $i = 1,2,3$. For each of these $i$, it turns out the Gr\"obner bases of $Az \cup P\cup f_i$ is $\{1\}$; in particular, the set $R$ is resolving. This result is not surprising since we showed previously that no non-trivial solutions exist for this system.

\subsection{Reduced Gr\"obner Basis of $P$}

We can improve the efficiency of Gr\"obner basis computations by examining the highly-structured set of polynomials $P$ in equation~(\ref{def:P}). In fact, an in depth look at the structure of $P$ provides a powerful intuition about how the complexity of the basis changes with the parameters $k$ and $a$ for the underlying Hamming graph $H_{k,a}$.

We can partition $P$ into disjoint sets (which we call blocks) as follows. For each $i = 0,\ldots,(k-1)$, define:
\begin{equation}
P_i :=\left\{\begin{array}{c}
    z_{ia+1}(z_{ia+1}-1)(z_{ia+1}+1)\\
    z_{ia+2}(z_{ia+2}-1)(z_{ia+2}+1)\\
    \vdots\\
    z_{(i+1)a}(z_{(i+1)a}-1)(z_{(i+1)a}+1)\\
    \hline 
    \sum\limits_{j = ia+1}^{(i+1)a} z_j\\
    \hline 
    (\sum\limits_{j = ia+1}^{(i+1)a} z_j^2)(2-\sum\limits_{j = ia+1}^{(i+1)a} z_j^2)\\
\end{array}\right.
\label{def:Pi}
\end{equation}

Note that $P = (P_0\cup \ldots \cup P_{k-1})$. Furthermore, each $P_i$ is a set of polynomials in the variables $z_{ia+1},\ldots, z_{(i+1)a}$, i.e. the $i$-th block of $z = \big((z_1,\ldots,z_a),\ldots,(z_{(k-1)a+1},\ldots,z_{ka})\big)$. In particular, no variable and hence no polynomial is shared between the blocks of $P$. Our next result reveals the importance of this partition.

In what follows, all polynomial ideals are with respect to $z$, not only variables in certain blocks of $z$.\\

\begin{lemma}[\cite{Cox_Little_OShea:2015}]
Consider two polynomial sets $Q_1$ and $Q_2$ with disjoint variables $x$ and $y$, and reduced Gr\"obner bases $G_1 \neq \{1\}$ and $G_2 \neq \{1\}$, respectively. Then $(G_1\cup G_2)$ is a reduced Gr\"obner basis of $(Q_1\cup Q_2)$.
\label{lem:groeb_block}
\end{lemma}

\noindent \emph{Proof: } The core part of this lemma is the disjoint variable vectors $x$ and $y$. Define $z:=(x,y)$, $Q:=(Q_1\cup Q_2)$, and $G:=(G_1\cup G_2)$.

We first show that $G$ is a Gr\"obner bases for $Q$. To start, note that $I(G)=I(Q)$. Indeed, since the ideal of a union is the sum of the ideals, and $I(Q_i)=I(G_i)$ for $i=1,2$ because $G_i$ is a Gr\"obner bases for $Q_i$, it follows that $I(Q)=I(Q_1)+I(Q_2)=I(G_1)+I(G_2)=I(G)$.

Next, we show that if $p,q\in I(G)$ then $Spoly(p,q)\stackrel{G}{\rightarrow} 0$. This is certainly the case if $p,q\in G_1$ or $p,q\in G_2$ because $G_1$ and $G_2$ are Gr\"obner bases. Thus, without loss of generality, we may assume that $p\in G_1$ and $q\in G_2$. In particular, $p=p(x)$ and $q=q(y)$, hence $LCM(LM(p),LM(q)) = LM(p)LM(q)$. As a result:
\begin{align*}
    Spoly(p,q) &= \frac{LM(p)LM(q)}{LT(p)}p - \frac{LM(p)LM(q)}{LT(q)}q \\
    &= \frac{LM(q)}{LC(p)}p - \frac{LM(p)}{LC(q)}q
\end{align*}
Then we write our polynomials as:
\begin{align*}
    p &= LT(p) + f \\
    q &= LT(q) + g 
\end{align*}
We can then write $Spoly(p,q)$ as:
\begin{align*}
    Spoly(p,q) &= \frac{LM(q)}{LC(p)}p - \frac{LM(p)}{LC(q)}q \\
               &= \frac{LM(q)f}{LC(p)} - \frac{LM(p)g}{LC(q)} \\
               &= \frac{LT(q)f}{LC(q)LC(p)} - \frac{LT(p)g}{LC(q)LC(p)} \\
               &= \frac{(q-g)f - (p-f)g}{LC(q)LC(p)} \\
               &= \frac{qf -pg}{LC(q)LC(p)}
\end{align*}
This formulation makes the reduction clear. $qf$ is trivially reducible to $0$ by $q$ and $pg$ is trivially reducible to $0$ by $p$, therefore $Spoly(p,q)$ reduces to $0$.

It remains to show that $G$ is also reduced. For this note that $LC(g)=1$ for each $g\in G$ because $G$ is the union of reduced Gr\"obner bases. 

Finally, we show that if $p,q\in G$ are distinct then no monomial of $q$ is divisible by the $LT(p)$. This is clearly the case if $p,q\in G_1$ or $p,q\in G_2$ because $G_1$ and $G_2$ are both reduced Gr\"obner bases. To complete the proof, without loss of generality assume that $p\in G_1$ and $q\in G_2$. Note that $p,q \neq 1$, otherwise $G_1 = \{1\}$ or $G_2 = \{1\}$ by Hilbert's Weak Nullstellensatz (Theorem~\ref{thm:weak_null}). Then, since $p=p(x)$ and $q=q(y)$, the leading term in $p$ cannot divide any monomial of $q$.  \hfill$\square$ 

The last lemma easily extends by induction to any finite union of polynomial sets which do not share variables. This leads immediately to the following result.

In what remains of this section, $G_i$ denotes the reduced Gr\"obner bases of $P_i$ in equation~(\ref{def:Pi}).\\

\begin{corollary}
$G = (G_0\cup \ldots \cup G_{k-1})$ is the reduced Gr\"obner bases of $P$ in equation~(\ref{def:P}).
\label{cor:groeb_block}
\end{corollary}

Finally, note that each $P_i$ is---up to the change of variables $z_j \longrightarrow z_{ia+j}$ for $j = 1,\ldots,a$---equal to $P_0$. In particular, the reduced Gr\"obner bases $G_i$ of $P_i$ may be obtained simply by renaming the variables in $G_0$. This drastically reduces the computation time of Gr\"obner bases to address the resolvability of vertices in Hamming graphs. 

We can further recognize that for $H_{k,a}$, $P_0$ only changes with $a$. A change in $k$ corresponds to a change in the number of blocks whereas a change in $a$ corresponds to a change in the structure of the blocks themselves. This means we only need to study how $P_0$ and $G_0$ change with $a$ without worrying about different values of $k$. In this regard, we have the following result.\\


\begin{lemma}
The reduced G\"obner bases of $P_0$ with respect to the lexicographic ordering is
\[G_0= \left\{\begin{array}{c}
         \sum\limits_{i=1}^a z_i\\
         z_2(z_2-1)(z_2+1)\\
        \vdots\\
         z_{a}(z_{a}-1)(z_{a}+1)\\
        \hline
        z_2z_3(z_2+z_3)\\
        \vdots\\
        z_{a-1}z_a(z_{a-1}+z_a)\\
        \hline
        z_2 z_3 z_4\\
        \vdots\\
        z_{a-2}z_{a-1}z_a
    \end{array}\right.\]
\end{lemma}

\noindent \emph{Proof:} 
$G_0$ can be shown to be a reduced Gr\"obner basis for an ideal by checking the 4 properties of reduced Gr\"obner bases and verifying that $I(G_0) = I(P_0)$. Demonstrating that $G_0$ satisfies the 4 properties of reduced Gr\"obner bases is a lengthy and uninformative process so we instead only show $I(G_0) = I(P_0)$.

To prove $I(G_0) = I(P_0)$ we show that the set of solutions to $G_0$ is the same as for $P_0$. The first polynomial simply requires that the sum of the variables is $0$. The next set of polynomials of the form $z_i(z_i-1)(z_i+1) = 0$ for $i = 2,3,\ldots,a$ forces $z_i \in \{-1,0,1\}$. The second set of polynomials $z_i^2 z_j + z_i z_j^2$ are formed for all pairs $z_i, z_j$ when $2\leq i<j$. This block enforces that for every pair of variables $z_i,z_j$, either $z_i = -z_j$ or at least one of the variables is $0$. The last set of polynomials $z_i z_j z_k$ are formed over all triplets $z_i,z_j,z_k$ when $2 \leq i < j < k$. This set of polynomials enforce that for every triplet of variables, at least one of them is $0$. Since every triplet must contain at least one $0$ element, there can be at most only one pair of non-zero elements. The constraints implied by these blocks are therefore identical to the constraints from Theorem~\ref{thm:Gen_System} and $I(G_0) = I(P_0)$. \hfill$\Box$

Another optimization can be made during the computation of the Gr\"obner bases $G_i$ of $G \cup f_i$. This basis computation can be improved by first performing the reduction $f_i \xrightarrow{G} r_i$ to speed up the subsequent S-polynomial computations and reductions. The polynomials $f_i$ are all very similar. In particular, if $f = \sum\limits_{j = 1}^{ak} z_j^2$ then $f_i = f - 2i$ for $i = 1,2,\ldots,k$. Using the following lemma, we can efficiently perform the reduction $f_i \xrightarrow{G} r_i$ for all $f_i$ with just one reduction $f \xrightarrow{G} r$.\\

\begin{lemma}
Consider a reduced Gr\"obner basis $G \neq \{1\}$ and a polynomial $f$ with no constant monomials. Define $f_i = f-c_i$ with all $c_i$ constants. If $f \xrightarrow{G} r$ then $f_i \xrightarrow{G} r-c_i$.
\end{lemma}

\noindent \emph{Proof:} Start by recognizing that $G \neq \{1\}$ implies $1 \notin I(G)$ and therefore $c_i \notin I(G)$ for any constant $c_i$. This implies that $c_i$ is irreducible by $G$, i.e. $c_i \xrightarrow{G} c_i$. 

If $f \xrightarrow{G} r$, we can write $f$ in terms of its quotients $q_j$.
\begin{align*}
    f &= q_1 g_1 + q_2 g_2 + \ldots + q_n g_n + r \\
\end{align*}
We also know that $f_i = f - c_i$ so we have:
\begin{align*}
    f_i &= f - c_i \\
        &= q_1 g_1 + q_2 g_2 + \ldots +  q_n g_n + r - c_i
\end{align*}
Since both $c_i$ and $r$ are irreducible by $G$, we know that $r-c_i$ is also irreducibly by $G$. This implies that $f_i = (q_1 g_1 + q_2 g_2 + \ldots q_n g_n + r - c_i)$ is a valid reduction of $f_i$ by $G$. This reduction is unique since $G$ is a Gr\"obner basis and so $f_i \xrightarrow{G} (r-c_i)$. $\square$ 

Using this lemma we arrive at the following result:
\begin{corollary}
Let $G$ be the reduced Gr\"obner basis of $P$ from equation~(\ref{def:P}) and $f = \sum\limits_{j = 1}^{ak} z_j^2$ with $f \xrightarrow{G} r$. If $f_i := (f-2i)$ for $i = 1,2,\ldots,k$ then $f_i \xrightarrow{G} (r-2i)$.
\end{corollary}

\section{Simplification on the Hypercubes}

In this section we specialize the results from the previous sections to hypercubes, i.e. Hamming graphs with an alphabet of size $a=2$. In this case, our findings and complexity simplify considerably to characterize resolving sets in $H_{k,2}$. We note that this simplification reproduces the result in~\cite{Beardon:2013}.

Our next result is a simplified version of Theorem~\ref{thm:Gen_System}.\\

\begin{corollary}
Let $R=\{v_1,\ldots,v_n\}$ be a set of nodes in $H_{k,2}$, and $\bar{v_i}$ denote the logical negation (flip) of $v_i$. Consider the matrix of dimensions $(n\times k)$ defined as:
\[B := \left(\begin{array}{c} 
v_1-\bar{v_1} \\ \vdots \\ v_n-\bar{v_n} 
\end{array}\right).\]
Then, $R$ resolves $H_{k,2}$ if and only if the equation $Bz=0$, with $z\in\{-1,0,1\}^k$, has only a trivial solution.
\label{cor:mat_system}
\end{corollary}

\noindent \emph{Proof:}
Recall the polynomial system $P=0$ from Section 2.1, specialized to the hypercube of dimension $k$:
\[\left\{\begin{array}{rcl}
    z_1(z_1-1)(z_1+1) &= 0 \\
    z_2(z_2-1)(z_2+1) &= 0 \\
    \vdots & \\
    z_{2k-1}(z_{2k-1}-1)(z_{2k-1}+1) &= 0 \\
    z_{2k}(z_{2k}-1)(z_{2k}+1) &= 0 \\
    \hline
    z_1+z_2 &= 0 \\
    \vdots & \\
    z_{2k-1}+z_{2k} &= 0 \\
    \hline
    (z_1^2+z_2^2)(2-(z_1^2+z_2^2)) &= 0 \\
    \vdots & \\
    (z_{2k-1}^2+z_{2k}^2)(2-(z_{2k-1}^2+z_{2k}^2)) &= 0
\end{array}\right.\]
Notice that $z_1+z_2 = 0$ implies that $z_1 = -z_2$ and $z_1^2 = z_2^2$; in particular, $z_1(z_1-1)(z_1+1) = -z_2(z_2-1)(z_2+1)$. So the equation $z_1(z_1-1)(z_1+1)=0$ is redundant and can be removed from the above polynomial system. Similarly, we find that
$(z_1^2+z_2^2)(2-(z_1^2+z_2^2))=-4z_2\cdot z_2(z_2-1)(z_2+1)$, in particular, the equation $(z_1^2+z_2^2)(2-(z_1^2+z_2^2))=0$ is also redundant and can be removed from the system too. The same reasoning can be followed starting with any of the equations in the middle block of the system. As a result, the polynomial system is equivalent to:
\[\left\{\begin{array}{rcl}
    z_2(z_2-1)(z_2+1) &= 0 \\
    z_1 + z_2 &= 0\\
    \vdots & \\
    z_{2k}(z_{2k}-1)(z_{2k}+1) &= 0 \\
    z_{2k-1}+z_{2k} &= 0
\end{array}\right.\]

But notice that all equations of the form $z_{2i-1} + z_{2i} = 0$ with $i=1,\ldots,k$ are linear constraints, so they can be pulled out of the polynomial system and instead encoded within the linear system $Az = 0$ in Theorem~\ref{thm:Gen_System}. This gives the simplified polynomial system:
\[\left\{\begin{array}{rcl}
    z_2(z_2-1)(z_2+1) &= 0 \\
    \vdots & \\
    z_{2k}(z_{2k}-1)(z_{2k}+1) &= 0 \\
\end{array}\right.\]
Note that this system is equivalent to imposing that $(z_2,\ldots,z_{2k})\in\{-1,0,1\}^k$.

Finally, we can encode the linear constraints $z_{2i-1} + z_{2i} = 0$ with $i=1,\ldots,k$ by setting $z = (-z_2, z_2,\ldots, -z_{2k}, z_{2k})^T$. Now consider the matrix $A$ from Theorem~\ref{thm:Gen_System} whose $i$-th column is denoted $A_i$:
\[A := \left(\begin{array}{ccccc}
\vline & \vline &  & \vline & \vline \\
A_1 & A_2 & \ldots & A_{2k-1} & A_{2k} \\
\vline & \vline & & \vline & \vline
\end{array}\right).\]
Then $Az = z_2(A_2-A_1) +\ldots + z_{2k}(A_{2k}-A_{2k-1})$, hence the linear system $Az=0$ can be represented by the alternative linear system: 
\[\begin{bmatrix}
        \vline &  & \vline \\
        (A_2 - A_1) & \ldots & (A_{2k}-A_{2k-1}) \\
        \vline & & \vline \\
    \end{bmatrix}
    \begin{bmatrix}
    z_2 \\
    \vdots \\
    z_{2k}
    \end{bmatrix}=0.\]
Recall however that from Theorem~\ref{thm:Gen_System}, the rows of $A$ are vectorized one-hot encodings; in particular, every pair of columns $A_{2i},A_{2i-1}$ corresponds to the $i$-th positions of the strings in the resolving set. This implies that $A_{2i} = v_i$ and $A_{2i-1} = \bar{v_i}$, and the corollary follows. \hfill$\square$
\newpage
\section{Applications of Results}

We are now ready to tackle the motivating question of this work: 
\begin{center}
\textit{Given a resolving set $R$ on $H_{k,a}$, what vertices, if any, can be removed such that $R$ remains resolving?} 
\end{center}
We can determine whether a vertex $r \in R$ can be removed by checking if the set $R\setminus r$ is resolving. If it is, then $r$ is removed, and this process is repeated until no nodes can be removed from $R$. Algorithm~\ref{alg:reduce} shows an improved version of this idea. Since $H_{k,a}$ is fixed, the polynomial set $P$ and the polynomials $f_i$ are also fixed. Therefore, it is more efficient to pre-compute the reduced Gr\"obner bases $G_i$ of $P \cup f_i$ before including the equations associated with the linear system $Az=0$. Additionally, each row of $A$ corresponds to a node in $R$. In particular, if $rank(A) < |R|$ then there are some linearly dependent rows of $A$. These linearly dependent rows correspond to vertices that can be removed from $R$ without changing the system. Since the system does not change upon their removal, these can be immediately removed beforehand with no checks of resolvability.

\begin{algorithm}
\begin{small}
\caption{Remove Vertices from Resolving Set}
\label{reduce_pseudo}
    \begin{algorithmic}
        \Function{reduceResolvingSet}{R,$k$,$a$}
            \State Construct polynomial set P and polynomials $f_i$ for $H_{k,a}$
            \State Construct matrix $A$ for set R
            \If{rank(A) < |R|}
                \State remove linearly dependent rows of A and corresponding vertices in R
            \EndIf
            \State Initialize shuffledR to be a randomly shuffled R
            \State $G_p$ $=$ precomputed reduced Gr\"obner basis of P
            \For{i = 1,2,\ldots,k}
                \State $G_i$ $=$ reduced incremental Gr\"obner basis of $G_P \cup f_i$
            \EndFor
            \For{r in shuffledR}
                \State remove r from R
                \State remove row corresponding to r from A
                \State $L$ = linear equations associated with $\hbox{rref}(A)z$
                \For{i = 1,2,\ldots,k}
                    \State $G$ $=$ $G_i$
                    \State $G$ $=$ reduced incremental Gr\"obner basis of $G \cup L$
                    \If{$G$ != \{1\}}
                        \State Add removed row back into A and r into R
                        \State \textbf{Break}
                    \EndIf
                \EndFor
            \EndFor
            \State \textbf{Return} R
        \EndFunction
    \end{algorithmic}
    \label{alg:reduce}
    \end{small}
\end{algorithm}

Unfortunately, the top-down approach of Algorithm~\ref{alg:reduce} may be very time consuming because removing nodes amounts to removing polynomials from the Gr\"obner bases $G_i$. Removing polynomials from a Gr\"obner basis is significantly more challenging than adding a polynomial to it. With this insight, it is more efficient to instead consider the bottom-up approach of randomly adding nodes from $R$ to an initially empty set $R'$, until $R'$ becomes resolving. The pseudocode for this approach can be found in Algorithm~\ref{alg:gen_reduce}.

\begin{algorithm}
\begin{small}
\caption{Generative Approach for Reducing Resolving Sets}
\label{reduce_pseudo}
    \begin{algorithmic}
        \Function{ReduceResolvingSet-Gen}{R,$k$,$a$}
            \State Construct polynomial system P and polynomials $f_i$ for $H_{k,a}$
            \State Construct matrix $A$ for set R
            \If{rank(A) < |R|}
                \State remove linearly dependent rows of A and corresponding vertices in R
            \EndIf
            \State Initialize shuffledR to be a randomly shuffled R
            \State $G_p$ $=$ precomputed reduced Gr\"obner basis of P
            \For{i = 1,2,\ldots,k}
                \State $G_i$ $=$ reduced incremental Gr\"obner basis of $G_P \cup f_i$
            \EndFor
            \State Initialize R' = $\emptyset$
            \For{r in shuffledR}
                \State add r to R'
                \State Set L = linear equation associated to r
                \For{i = 1,2,\ldots,k}
                    \State $G_i$ $=$ $G_i \cup L$
                \EndFor
                \If{All $G_i == \{1\}$}
                    \State \textbf{Return} R'
                \EndIf
            \EndFor
        \EndFunction
    \end{algorithmic}
    \label{alg:gen_reduce}
    \end{small}
\end{algorithm}

The insights for reducing known resolving sets can be used to generate resolving sets from scratch as well. Algorithm~\ref{alg:gen} implements this approach. The algorithm incrementally adds random nodes to an initially empty set $R$ that increase the rank of the associated matrix $A$. This process is repeated until $R$ is resolving. The algorithm guarantees construction of a resolving set since $A$ is non-singular (i.e., has only the trivial solution) when $rank(A) = ak$. Since $rank(A) = |R|$, once $ak$ vertices have been added into $R$, it is trivially resolving on $H_{k,a}$. We note that the resolving sets produced by this algorithm are not necessarily minimal.

\begin{algorithm}
\begin{small}
\caption{Generate Resolving Sets}
\label{GenR_pseudo}
    \begin{algorithmic}
        \Function{genResolvingSet}{$k$,$a$}
            \State Construct polynomial system P and polynomials $f_i$ for $H_{k,a}$
            \State Precompute reduced Gr\"obner bases $G_i$ for $P \cup f_i$
            \State Initialize $V$ $:=$ vertex set of $H_{k,a}$
            \State $R$ = some randomly selected $r \in V$
            \While{R is not resolving}
                \State Randomly select some $r \in V$ that has not been selected before
                \State Construct matrix $A$ for $R \cup r$
                \If{rank(A) == |R|+1}
                    \State $R$ = $R \cup r$
                    \For{i = 1,2,\ldots,k}
                        \State $G$ = reduced incremental Gr\"obner basis of $G_i \cup \hbox{rref}(A)z$
                        \If{G != \{1\}}
                            \State \textbf{Break}
                        \EndIf
                        \State Set $R$ as resolving
                    \EndFor
                \EndIf
            \EndWhile
            \State \textbf{Return} R
        \EndFunction
    \end{algorithmic}
    \label{alg:gen}
    \end{small}
\end{algorithm}

\section{Preliminary run-time Analysis}

We implemented the proposed Gr\"obner basis approach for checking resolvability on $H_{k,a}$ using SymPy~\cite{SymPy}, a Python computer algebra package. In this section we present run-time results of this implementation. 

To test the theory, we have generated $4,359$ example sets for Hamming graphs with parameters $k = 1,\ldots,10$ and $a = 2,\ldots,5$. The graphs were restricted to those where $ak \leq 25$ in order to limit the complexity for these initial tests. The example sets consist of around 200 sets per graph, at least half of which are resolving on their associated Hamming graph. 

\begin{figure}[h]
    \centering
    \includegraphics[width = 0.6\textwidth]{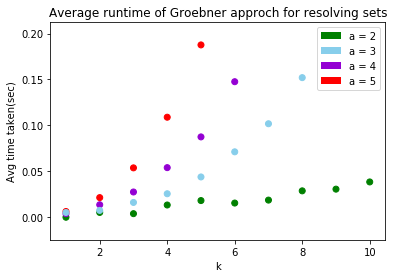}
    \caption{run-time of checking resolvability for resolving sets.}
    \label{fig:groeb_resolving}
\end{figure}

Figure~\ref{fig:groeb_resolving} displays the average run-time for showing a set is resolving in $H_{k,a}$ (for various combinations of $k$ and $a$) using the developed Gr\"obner basis algorithm. In this plot, there is a clear and promising polynomial growth in $k$, with almost linear behavior for $a = 2$. Additionally, the actual run-times are extremely impressive. All times are below $0.2$ seconds---even for $H_{8,3}$ which contains over $6,500$ vertices. Note however that this is specifically for sets which are resolving on their associated Hamming graphs. As we will see throughout this section, resolving sets seem to produce linear systems with simple reduced Gr\"obner bases, drastically improving both the average run-time as well as the time consistency of our algorithm.

\begin{figure}[h]
    \begin{subfigure}[b]{0.48\textwidth}
        \centering
        \includegraphics[width = 0.99\textwidth]{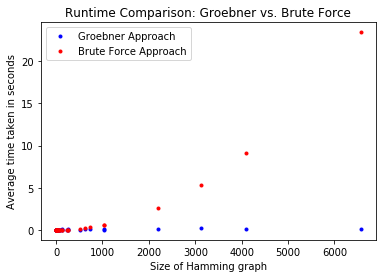}
        \caption{Gr\"obner basis algorithm significantly outperforms brute-force approach for resolving sets.}
    \end{subfigure}
    \hfill
    \begin{subfigure}[b]{0.48\textwidth}
        \centering
        \includegraphics[width = 0.99\textwidth]{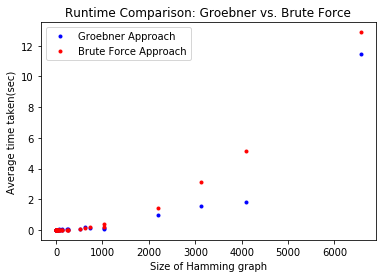}
        \caption{Gr\"obner basis algorithm performs comparably to brute-force when also considering non-resolving sets}
    \end{subfigure}
    \caption{}
    \label{fig:comparison}
\end{figure}

Next, we consider how the run-time of the brute-force algorithm and Gr\"obner basis algorithm compare as the size of the Hamming graph increases. As seen in Figure~\ref{fig:comparison}, the Gr\"obner basis algorithm vastly outperforms the brute-force algorithm for showing a set is resolving, boasting over $100$ times faster results. By contrast, the Gr\"obner basis algorithm performs comparably to the brute-force algorithm for general sets where some are non-resolving. This shows that the Gr\"obner basis algorithm is faster for resolving sets than for non-resolving sets whereas the opposite is true with respect to the brute-force algorithm. A possible explanation is that resolving sets are structurally simpler than non-resolving sets. This results in a simpler system $P \cup Az = 0$ for resolving sets and the subsequent Gr\"obner basis computations are very quick. Non-resolving sets instead produce systems which have more than just the trivial solution whose associated Gr\"obner bases are much more cumbersome and complex to handle. The brute-force algorithm by comparison is much faster for non-resolving sets since it works by finding an unresolved pair of vertices. If no unresolved pair of vertices exist, the brute-force algorithm must go through all possible pairs, increasing the run-time. 

\begin{figure}[h]
    \centering
    \includegraphics[width = 0.6\textwidth]{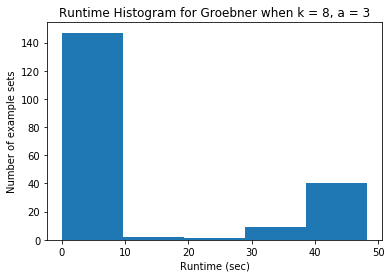}
    \caption{Histogram of Gr\"obner basis run-times on $H_{8,3}$. Almost $\frac{3}{4}$-ths of the trials took under $10$ seconds with the last $\frac{1}{4}$ taking around $45$ seconds.}
    \label{fig:groeb_hist}
\end{figure}

Finally, we analyze the results of the Gr\"obner basis algorithm to get a better understanding of why the algorithm is slower on non-resolving sets. In Figure~\ref{fig:groeb_hist} there is a striking range of run-time results. While the vast majority of trials took under $10$ seconds, one quarter of example sets took around $45$ seconds. It turns out that those example sets were not resolving and were constructed with similar nodes, yielding similar linear systems and solutions. There are two possible explanations for this stark contrast in performance. As before, the linear systems constructed from these resolving sets may produce highly complex solution spaces which makes Gr\"obner basis computations very slow. It is odd however that there are almost no trials that land in between these two peaks in run-time. Another possible explanation is that SymPy Gr\"obner basis computations are unstable on these systems. Indeed, there are documented and unresolved issues with SymPy Gr\"obner basis computations where two systems of similar complexity take vastly different times to complete. The quarter of example sets with significantly higher runtime are very similar and produce nearly identical linear systems. It is possible that all of these systems were unstable for SymPy in the same ways. This could also explain the empty gap of run-times since systems are either drastically unstable and take a very long time to compute or they are stable and are completed within $10$ seconds. A combination of both effects could certainly be at play, motivating potential future study to isolate the two different effects.

\section{Discussion and Future Work}

We have shown that resolvability on Hamming graphs can be characterized in terms of the solution space of a polynomial system specific to the Hamming graph, and a linear system constructed from the elements of a set of nodes $R$ (Theorem~\ref{thm:Gen_System}). This characterization provides a novel framework with which to analyze resolving sets, connecting the theory of metric dimension on Hamming graphs to algebraic geometry and linear algebra. Our formulation can be simplified dramatically in hypercubes (Corollary~\ref{cor:mat_system}). Through this simplification, we are able to retrieve the same linear system as Beardon~\cite{Beardon:2013}, suggesting that our polynomial system somehow contains all the logical constraints implied by a binary alphabet. The characterization of resolvability has immediate utility in improving resolving set based applications such as k-mer embeddings, but it also illuminates a deep connection between algebraic geometry and Hamming graphs.

We have also exploited the block structure and symmetry of the polynomial system we encountered to explicitly define the Gr\"obner basis for any Hamming graph (Corollary~\ref{cor:groeb_block}). Further exploring the block structure of the Gr\"obner bases could reveal potentially powerful optimizations of incremental basis computations, or further insights into resolvability at a structural level. Since these bases are nearly identically structured for all Hamming graphs, it might be possible to start with the Gr\"obner basis and try to derive a linear system which would cause it to reduce to $\{1\}$. Discovering such a derivation would amount to an explicit construction of resolving sets and possibly be the basis of proving the minimality of resolving sets. 

Finally, we have analyzed the performance of the Gr\"obner basis theory for checking resolvability. The new method drastically outperforms the brute-force approach when the sets are resolving but there are significant reductions in speed for non-resolving sets. The specifics of these slow-downs suggest either highly complex solution spaces to the constructed linear systems or instability in the SymPy Gr\"obner basis computations for these specific examples. Testing this implementation in other computer algebra systems will be an important step in determining the root cause of the severe decrease in performance. Even with the potential SymPy instability, our new Gr\"obner basis algorithm performs better than brute-force on average while also providing a structural insight into resolvability on Hamming graphs.
\bibliographystyle{ieeetr}
\bibliography{main} 
\end{document}